\documentclass[12pt,oneside,english]{amsart}
\usepackage[latin1]{inputenc}
\usepackage{amssymb}

\makeatletter

\newtheorem{lemma}{Lemma}

\newtheorem{definition}{Definition}

\newtheorem{remark}{Remark}

\newtheorem{problem}{Problem}

\numberwithin{equation}{section}
\usepackage{babel}

\makeatother
\begin{document}
\baselineskip=17pt

\title[The m\'{e}nage problem with a fixed couple]{The m\'{e}nage problem with a
fixed couple}

\author{Vladimir Shevelev}
\address{Department of Mathematics \\Ben-Gurion University of the
 Negev\\Beer-Sheva 84105, Israel. e-mail:shevelev@bgu.ac.il}
\author{Peter J. C. Moses}
 \address{United Kingdom. e-mail: mows@mopar.freeserve.co.uk}


\begin{abstract}
We give a solution of the following variation of the classic m\'{e}nage problem:
"Let one from $n$ married couples in the m\'{e}nage problem (see Problem 1) be a
couple of $M$ and his wife. After the ladies are seated at every other chair,
$M$ (in token of respect) is the first man allowed to choose one of the remaining
chairs. To find the number of ways of seating the other men, with no man seated next
to his wife, if $M$ chooses the chair that is $d$ seats clockwise from his wife's
chair."
 \end{abstract}

\maketitle

\section{Introduction }
  In 1891, Lucas \cite{2} formulated the following "m\'{e}nage problem":
\begin{problem}\label{pr1}
   To find the number $M_n$ of ways of seating $n$ married couples at a circular table, men and women in alternate positions, so that no husband is next to his wife.
\end{problem}
    After seating the ladies by $2n!$ ways we have
 \begin{equation}\label{1.1}
 M_n=2n!U_n,
\end{equation}
where $U_n$ is the number of ways of seating men.\newline
\indent Earlier Muir \cite{4} solved a problem posed by Tait (cf. \cite{4}): to find the number $H_n$ of
 permutations $\pi$ of $\{1,...,n\}$ for which $\pi(i)\neq i$ and $\pi(i)\neq i+1 \pmod {n}, \enskip i=1,...,n.$  Simplifying Muir's solution, Cayley \cite{1} found a very simple recursion for $H_n: H_2=0, H_3=1,$ and for $n\geq4,$
 \begin{equation}\label{1.2}
 (n-2)H_n=n(n-2)H_{n-1}+nH_{n-2}+4(-1)^{n+1}.
\end{equation}
 13 years later, Lucas \cite{2} on page 495, formula (9), gave the same formula
 for $U_n.$  So we see that \cite{1}-\cite{2} imply the equality
\begin{equation}\label{1.3}
 H_n=U_n
\end{equation}
which became well known after development of the rook theory \cite{5}, ch. 7,8.
In Section 2 we give also a one-to-one correspondence between $H_n$ and $U_n.$
In 1934, Touchard \cite{8} found a remarkable explicit formula
\begin{equation}\label{1.4}
U_n=\sum_{k=0}^{n}(-1)^k \frac{2n} {2n-k}\binom {2n-k} {k}(n-k)!
\end{equation}
A beautiful proof of (\ref{1.4}) with help of the "rook technique" one can find in \cite{5}.
The first terms of the sequence $\{U_n\},$ for $n\geq2,$ are (cf. A000179 in
\cite{7})
 \begin{equation}\label{1.5}
0, 1, 2, 13, 80, 579, 4738, 43387, 439792, 4890741, 59216642, ...
\end{equation}
Note that formulas for $U_n$ in other forms are given by Wayman and Moser \cite{9}
and Shevelev \cite{6}.\newline
In the present paper we study the following problem.
\begin{problem}\label{pr2}
Let one from n married couples in the m\'{e}nage problem (see Problem \ref{pr1}) be
a couple of $M$ and his wife. After the ladies are seated at every other chair,
$M$ is the first man allowed to choose one of the remaining chairs. To find the
number of ways of seating the other men, with no man seated next to his wife, if $M$
chooses the chair that is $d$ seats clockwise from his wife's chair.
\end{problem}

\section{Equivalence of Tait's and m\'{e}nage problems}
Let $A$ be $n\times n$ $(0,1)$-matrix. To every permutation $\pi=\{i_1,...,i_n\}$ of numbers $\{1,...,n\}$ corresponds a set of positions of $A$ $\{(1,i_1),...,(n,i_n)\},$ which is called a diagonal of $A.$ Thus matrix $A$ has $n!$ distinct diagonals. If to place in every position a chessboard piece rook, then we have $n$ non-taking rooks. The number of distinct ways of putting $n$ non-taking rooks on \slshape positions of 1's \upshape \enskip of matrix $A$ is called \slshape permanent \upshape \enskip of $A$ ($per A$). Denote $J=J_n$ $n\times n$ matrix which consists of 1's only. It is clear that $per J_n=n!$ Denote, furthermore, $I=I_n$ and $P=P_n$ $n\times n$ $(0,1)$-matrices $(n\geq3)$ each of which contains only diagonal of 1's: $(1,1),...,(n,n)$ and $(1,2),...,(n-1,n), (n,1)$ correspondingly. In Tait's problem we should find the number of permutation with the prohibited positions $(1,1),...,(n,n)$ and $(1,2),...,(n-1,n), (n,1).$ Therefore, Tait's problem is the problem of calculation of $H_n=per(J_n-I-P).$\newline
\indent Consider now the m\'{e}nage problem. Denote $2n$ chairs at a circular table
 by the symbols
\begin{equation}\label{2.1}
1,\overline{1}, 2, \overline{2},..., n,\overline{n}
\end{equation}
clockwise. Ladies occupy either chairs $\{1,...,n\}$ or chairs $\{\overline{1},...,\overline{n}\}.$ Let they occupy, say, chairs $\{\overline{1},...,\overline{n}\}.$ Then to every man we give a number $i,$
 if his wife occupies the chair $\overline{i}.$ Now the $i$-th man,
for $i=1,...,n-1,$ can occupy every chair except of chairs $i, i+1,$ while the
$n$-th man cannot occupy chairs $n$ and $1.$ Denoting in the corresponding
$n\times n$ incidence matrix the
\newpage prohibited positions by 0's and other positions
by 1's, we again obtain the matrix $J_n-I-P.$  Now, evidently, to every seating
the men corresponds a diagonal of 1's in this matrix. This means that
(cf. \cite{3})
\begin{equation}\label{2.2}
 U_n=per(J_n-I-P),
\end{equation}
and (\ref{1.3}) follows.\newline
\indent Moreover, we can indicate a one-to-one correspondence between the diagonals
of $1's$ of the matrix $J_n-I-P$ and arrangements of $n$ married couples around a
circular table by the rules of the m\'{e}nage problem, after the ladies
$w_1, w_2, ..., w_n$ have taken the chairs numbered
\begin{equation}\label{2.3}
2n, 2, 4, ..., 2n-2
\end{equation}
respectively. Suppose we consider a diagonal of $1's$ of the matrix $J_n-I-P$:
\begin{equation}\label{2.4}
(1, j_1), (2, j_2), ..., (n,j_n)).
\end{equation}
Then the men $m_1, m_2, ..., m_n$ took chairs with numbers
\begin{equation}\label{2.5}
2j_i-3 \pmod {2n},\enskip i=1,2,...,n,
\end{equation}
where the residues are chosen from the interval $[1,2n].$\newline
Indeed $\{j_i\}$ is a permutation of $1,...,n.$ So $\{2j_i-3\} \pmod {2n}$ is a
permutation of odd positive integers not exceeding $2n-1$. Besides, the distance
between places of $m_i$ (\ref{2.5}) and $w_i$ (\ref{2.3}) cannot be 1. Indeed, the
equality $|2(j_i-i)-1| = 1 \pmod{2n}$ is possible if and only if either $j_i=i$
or $j_i=i+1 \pmod {n}$ that correspond to positions of $0's$ in matrix $J_n-I-P.$
\newline For example, in case of $n=5$ and $j_1=3, j_2=1, j_3=5, j_4=2, j_5=4$ in
(\ref{2.4}), by (\ref{2.3}) and (\ref{2.5}), the chairs $1,2,...,10$ are taken by
$$m_4, w_2, m_1, w_3, m_5, w_4, m_3, w_5, m_2, w_1,$$ respectively.

\section{Equivalent form of Problem 2}
 A one-to-one correspondence between the diagonals of $1's$ of the matrix $J_n-I-P$
and arrangements of $n$ married couples around a circular table by the rules of the
m\'{e}nage problem, after the ladies $w_1, w_2, ..., w_n$ have taken the chairs
numbered as in (\ref{2.3}) prompts us a way
of solution of Problem \ref{pr2}. Indeed, denote by $(J_n-I-P)[1|\enskip r]$ the
matrix obtained by the removing the first row and the $r$-th column of the matrix
$J_n-I-P.$  Then, by expansion of the permanent (\ref{2.2}) over the first row, we
have
\begin{equation}\label{3.1}
U_n=\sum_{r=3}^{n}per((J_n-I-P)[1|\enskip r] ).
\end{equation}
\newpage
Let $M=m_1.$ Since, in view of symmetry, the solution is invariant from a chair that
$m_1's$ wife $w_1$ occupy, let $w_1$ occupy the chair $2n$ (or $0 \pmod{2n}).$
Then, by (\ref{2.5}), in (\ref{3.1}) to $r=3,4,5,6,...$ there correspond the values
of distances
$d=3,5,7,9,...$  clockwise, i.e., the distance
\begin{equation}\label{3.2}
 d=2r-3, \enskip r\geq3,
\end{equation}
clockwise between $m_1$ and $w_1.$ Thus for the solution of Problem \ref{pr2} we
should find the summands of (\ref{3.1}). But technically this problem is rather
difficult. Here we can solve it only due to, as we show, a representation of rook
polynomials of each matrix $A_r=(J_n-I-P)[1|\enskip r], \enskip 3\leq r\leq n, $ as
a product of rook polynomials of simpler matrices.

\section{Lemmas}

\indent Let now $M$ be a rectangle (quadratic) $(0,1)$-matrix.
\begin{definition}\label{d1}
The polynomial
\begin{equation}\label{4.1}
 R_M(x)=\sum_{j=0}^{n}\nu_j(M)x^j
\end{equation}
where $\nu_0=1$ and $\nu_j$ is the number of ways of putting $j$ non-taking rooks on positions 1's of $M,$ is called \upshape rook polynomial.\slshape
\end{definition}
In particular, if $M$ is a quadratic $n\times n$-matrix, then $\nu_n(M)=per M.$\newline
\indent Now we formulate three main results (Lemmas \ref{L1}-\ref{L3}) of the classic Kaplansky-Riordan rook theory (cf. \cite{5}, Ch. 7-8).
\begin{lemma}\label{L1}
If $M$ is a quadratic matrix with the rook polynomial $(4.1),$ then
\begin{equation}\label{4.2}
per(J_n-M)=\sum_{j=0}^{n}(-1)^j\nu_j(M)(n-j)!
\end{equation}
\end{lemma}
\begin{definition}\label{d2}
Two submatrices $M_1$ and $M_2$ of $(0,1)$-matrix $M$ are called \upshape disjunct \slshape if no 1's of $M_1$ in the same row or column as those of $M_2.$
\end{definition}
From Definition \ref{d1} the following lemma evidently follows.
\begin{lemma}\label{L2}
If $(0,1)$-matrix $M$  consists of two disjunct submatrices $M_1$ and $M_2,$  then
\begin{equation}\label{4.3}
R_M(x)=R_{M_1}(x)R_{M_2}(x).
\end{equation}
\end{lemma}
Consider a position $(i,j)$ of 1 in matrix $M.$ Denote $M^{(0(i,j))}$ the matrix
obtained from $M$ after replacing 1 in position $(i,j)$ by 0. Denote $M^{(i,j)}$
the matrix obtained from $M$ by  removing the $i$-th row and $j$-column.
\newpage
\begin{lemma}\label{L3} We have
 \begin{equation}\label{4.4}
R_M(x)=xR_{M^{(i,j)}}+R_{M^{(0(i,j))}}.
\end{equation}
\end{lemma}

Consider so-called simplest connected staircase $(0,1)$-matrices. Such matrix is
called $k$-staircase, if the number of its 1's equals to $k.$
For example, the following several matrices are $5$-staircase:
$$\begin{pmatrix}1&1&0\\0&1&1\\0&0&1
\end{pmatrix}, \begin{pmatrix}1&0&0\\1&1&0\\0&1&1
\end{pmatrix}, \begin{pmatrix}0&0&0\\1&0&0\\1&1&0\\0&1&1
\end{pmatrix}, \begin{pmatrix}0&0&1&1&0\\0&1&1&0&0\\0&1&0&0&0
\end{pmatrix}$$
and the following matrices are $6$-staircase:
$$\begin{pmatrix}1&1&0&0\\0&1&1&0\\0&0&1&1
\end{pmatrix},\begin{pmatrix}1&0&0\\1&1&0\\0&1&1\\0&0&1
\end{pmatrix}, \begin{pmatrix}0&0&0\\1&0&0\\1&1&0\\0&1&1\\0&0&1
\end{pmatrix},\begin{pmatrix}0&0&1&1&0\\0&1&1&0&0\\1&1&0&0&0
\end{pmatrix}$$
\begin{lemma}\label{L4}
For every $k\geq1,$ all $k$-staircase matrices $M$ have the same rook polynomial
 \begin{equation}\label{4.5}
R_M(x)=\sum_{i=0}^{\lfloor\frac{k+1} {2}\rfloor}\binom {k-i+1}{i}x^i.
\end{equation}
\end{lemma}
\begin{proof}

For each $k$-staircase $(0,1)$-matrix, one can carry out the same review as for
the following simplest connected $k$-staircase $(0,1)$-matrix with the configuration
of 1's of the form:
$$\begin{matrix}1&1&.&.&.&.\\.&1&1&.&.&.\\.&.&.&.&.&.\\.&.&.&1&1&.\\.&.&.&.&1&1
\end{matrix}$$

The last right 1 is absent for odd $k$ and is present for even $k.$
In both cases, by Lemma \ref{L3}, for the rook polynomial $R_k(x)$
we have
$$R_0(x)=1,\enskip R_1(x)=x+1, \enskip R_k(x)=R_(k-1)(x) + xR_(k-2)(x),\enskip k>=2.$$
This equation has solution (\ref{4.5}) of Lemma \ref{L4} (see \cite{5}, ch.7,
eq.(27)).
\end{proof}

\section{Solution of Problem 2}
According to Lemma \ref{L1}, in order to calculate permanent of matrix
$(J_n-I-P)[1|\enskip r],$ we can find rook polynomial of matrix
$J_{n-1}-(J_n-I-P)[1|\enskip r].$ We use the equation
\begin{equation}\label{5.1}
J_{n-1}-(J_n-I-P)[1|\enskip r]=(I_n+P)[1|\enskip r].
\end{equation}
\newpage
which follows from the additivity $(A+B)[1|\enskip r]=A[1|\enskip r]+B[1|\enskip r]$
and equality $J_{n-1}=J_n[1|\enskip r].$

Pass from matrix $(I_n+P)$ to matrix $(I_n+P)[1|\enskip r].$  We have
(here $n=10,\enskip r=5$)
\begin{equation}\label{5.2}
\begin{pmatrix}1&1&0&0&0&0&0&0&0&0\\0&1&1&0&0&0&0&0&0&0\\0&0&1&1&0&0&0&0&0&0\\0&0&0&1&1&0&0&0&0&0
\\0&0&0&0&1&1&0&0&0&0\\0&0&0&0&0&1&1&0&0&0\\0&0&0&0&0&0&1&1&0&0\\0&0&0&0&0&0&0&1&1&0\\
0&0&0&0&0&0&0&0&1&1\\1&0&0&0&0&0&0&0&0&1
\end{pmatrix}\rightarrow \begin{pmatrix}0&1&1&0&0&0&0&0&0\\0&0&1&1&0&0&0&0&0\\0&0&0&1&0&0&0&0&0\\0&0&0&0&1&0&0&0&0
\\0&0&0&0&1&1&0&0&0\\0&0&0&0&0&1&1&0&0\\0&0&0&0&0&0&1&1&0\\0&0&0&0&0&0&0&1&1\\1&0&0&0&0&0&0&0&1
\end{pmatrix}
\end{equation}

Now we use Lemma \ref{L3} to the latter matrix in case $i=n-1,\enskip j=1.$ Denote
\begin{equation}\label{5.3}
A=((I_n+P)[1|\enskip r])^{(n-1,1)}, \enskip B=((I_n+P)[1|\enskip r])^{(0(n-1,1))}.
\end{equation}
According to (\ref{4.4}), we have
 \begin{equation}\label{5.4}
R_{(I_n+P)[1|\enskip r]}(x)=xR_{A}(x)+R_{B}(x).
\end{equation}
Note that matrix $A$ has the form (here $n=10, \enskip r=5$)
\begin{equation}\label{5.5}
A=\begin{pmatrix}1&1&0&0&0&0&0&0\\0&1&1&0&0&0&0&0\\0&0&1&0&0&0&0&0\\0&0&0&1&0&0&0&0
\\0&0&0&1&1&0&0&0\\0&0&0&0&1&1&0&0\\0&0&0&0&0&1&1&0\\0&0&0&0&0&0&1&1
\end{pmatrix}
\end{equation}
which is $(n-2)\times(n-2)$ matrix with $2n-6$ 1's.
This matrix consists of two disjunct matrices: $(r-2)\times(r-2)$ matrix $A_1$ of
the form (here $r=5$)
\begin{equation}\label{5.6}
A_1=\begin{pmatrix}1&1&0\\0&1&1\\0&0&1
\end{pmatrix}
\end{equation}
which is $2r-5$-staircase matrix, and $(n-r)\times(n-r)$ matrix
(here $n=10,\enskip r=5$)
\begin{equation}\label{5.7}
A_2=\begin{pmatrix}1&0&0&0&0\\1&1&0&0&0\\0&1&1&0&0\\0&0&1&1&0\\0&0&0&1&1
\end{pmatrix}
\end{equation}
\newpage
which is $2(n-r)-1$-staircase matrix.\newline
\indent Thus, by Lemmas \ref{L2} and \ref{L4}, we have

$$R_A(x)=\sum_{i=0}^{r-2}\binom {2r-i-4}{i}x^i\sum_{i=0}^{n-r}\binom {2(n-r)-i}{i}
x^i$$ \newline
\begin{equation}\label{5.8}
=\sum_{i=0}^{r-2}\binom {2r-i-4}{i}x^i\sum_{j=0}^{n-r+1}\binom {2(n-r)-j+1}{j-1}x^{j-1}.
\end{equation}

Note that, since $\binom {n}{-1}=0$, then we write formally the lower limit in
interior sum $j=0.$
Furthermore, $(n-1)\times (n-1)$ matrix $B=B(r)$ has the form (here $n=10, \enskip r=5$)

\begin{equation}\label{5.9}
B= \begin{pmatrix}0&1&1&0&0&0&0&0&0\\0&0&1&1&0&0&0&0&0\\0&0&0&1&0&0&0&0&0\\0&0&0&0&1&0&0&0&0
\\0&0&0&0&1&1&0&0&0\\0&0&0&0&0&1&1&0&0\\0&0&0&0&0&0&1&1&0\\0&0&0&0&0&0&0&1&1\\0&0&0&0&0&0&0&0&1
\end{pmatrix}
\end{equation}
and contains $2n-5$ 1's.
This matrix consists of two disjunct matrices: $(r-2)\times(r-1)$ matrix $B_1$ of
the form (here $r=5$)
\begin{equation}\label{5.10}
B_1=\begin{pmatrix}0&1&1&0\\0&0&1&1\\0&0&0&1
\end{pmatrix}
\end{equation}
which is $(2r-5)$-staircase matrix, and $(n-r+1)\times(n-r)$ matrix
(here $n=10,\enskip r=5$)
\begin{equation}\label{5.11}
B_2=\begin{pmatrix}1&0&0&0&0\\1&1&0&0&0\\0&1&1&0&0\\0&0&1&1&0\\0&0&0&1&1\\0&0&0&0&1
\end{pmatrix}
\end{equation}
which is $2(n-r)$-staircase matrix. Thus, by Lemmas \ref{L2} and \ref{L4}, we have \newline
\begin{equation}\label{5.12}
R_B(x)=\sum_{i=0}^{r-2}\binom {2r-i-4}{i}x^i\sum_{j=0}^{n-r+1}\binom {2(n-r)-j+1}{j}x^j.
\end{equation}\newpage
Note that, since $\binom {n-r}{n-r+1}=0$, then we write formally the upper limit
in interior sum $j=n-r+1.$
Now from (\ref{5.4}), (\ref{5.8}) and (\ref{5.12}) we find\newline
$$R_{(I_n+P)[1|\enskip r]}(x)=\sum_{i=0}^{r-2}\binom {2r-i-4}{i}x^i\sum_{j=0}^{n-r+1}\binom {2(n-r)-j+2}{j}x^j $$\newline
\begin{equation}\label{5.13}
=\sum_{k=0}^{n-1}x^k\sum_{i=0}^{\min (k,\enskip r-2)}\binom {2r-i-4}{i}\binom {2(n-r)-k+i+2}{k-i}.
\end{equation}\newline
Note that in the interior sum in (\ref{5.13}) it is sufficient to take summation
over interval $[\max(r+k-n-1,0), \min(k,r-2)].$  Thus, by Lemma \ref{L1} and
(5.1), we have

$$per((J_n-I-P)[1|\enskip r])=\sum_{k=0}^{n-1}(-1)^k (n-k-1)!\times$$\newline
\begin{equation}\label{5.14}
\sum_{i=\max(r+k-n-1,\enskip0)}^{\min (k,\enskip r-2)}\binom {2r-i-4}{i}\binom {2(n-r)-k+i+2}{k-i}.
\end{equation}\newline
Formula (\ref{5.14}) solves Problem \ref{pr2} for $r=(d+3)/2\geq3$ (by (3.2)) and naturally $n>(d+1)/2.$. $\blacksquare$\newline
The sequences corresponding to $d=3,5,7,9$ and 11 see in \cite{7}, A258664-A258667
and A258673.
\begin{remark}
The prohibited values $d=1$ and $d=2n-1$ correspond to the case when M and W are seated at neighboring chairs. Let us calculate the number of of ways of seating the remaining
men after that all the women occupied their chairs so that women and men are in alternate chairs but M and W are the only couple seated next to each other.
Thus we have a classic m\'{e}nage problem for $n-1$ couples for a one-sided linear
table, after that ladies already occupied their chairs. So, by \cite{5}, Ch. 8, Th. 1,
$t=0,$ the solution $V_n$ of this problem is
\begin{equation}\label{5.15}
V_n=\sum_{k=0}^{n-1}(-1)^k  \binom{2n-k-2}{k}(n-k-1)!, \enskip n>1.
\end{equation}
One can verify that this result is also formally obtained from $(\ref{5.14})$ for
both $r=2$ and $r=n+1.$  It could be proven easily independently. Cf. also $A259212$
in \cite{7}.
\end{remark}
\newpage 
\section{Acknowledgment}
The authors thank Giovanni Resta and Jon E. Schoenfield for useful discussions.

\;\;\;\;

\;\;\;\;\;\;\;\;

\end{document}